\documentclass{birkmult}
\usepackage[utf8]{inputenc}
\usepackage[T1]{fontenc}
\usepackage[russian,english]{babel}

\usepackage{url}

\title
{The zeros of random sections of real vector bundles}

\author{Boris Kazarnovskii}
\address{%
Institute for Information Transmission Problems of the Russian Academy of Sciences,
Moscow, Russia
}
\email{kazbori@gmail.com}

\newtheorem{theorem}{Theorem}
\theoremstyle{definition}
\newtheorem{definition}{Definiton}
\theoremstyle{lemma}
\newtheorem{lemma}{Lemma}
\theoremstyle{corollary}
\newtheorem{corollary}{Corollary}
\theoremstyle{definition}
\newtheorem{example}{Example}
\def\vol{{\rm vol}}
\def\R{{\mathbb R}}

\begin{document}
\begin{abstract}
We define integral geometric analogues of the Chern classes for
real vector bundle on a smooth real variety. More precisely, we
define the Chern densities of a real bundle. These densities are
analogues of the Chern forms of a complex vector bundle and inherit
some of their properties.

\noindent
(The text is a summary of a report on
the conference PCA'2024 in
Euler International Mathematical Institute, St. Petersburg)
\end{abstract}

\maketitle
\section*{The zeros of random systems of functions}
We begin with a theorem on the number of common zeros of random systems of functions from \cite{AK}.
Let $V$ be a finite-dimensional space of smooth functions on an $n$-dimensional differentiable manifold $X$.
Consider a random system of equations
\begin{equation}\label{eq1}
  f_1=\ldots=f_n=0,\,\, 0\ne f_i\in V
\end{equation}
Denote by $N(f_1,\ldots,f_n)$
the number of solutions of the system (\ref{eq1}).
We define the randomness of the system using a certain scalar product in $V$ as follows:
we consider the functions $f_i$ as independent random vectors in $V$
with respect to the Gaussian measure chosen in $V$ according to the chosen scalar product.
The situation with a more general choice of probability distribution in $V$ is described in \cite{Ka1}.
Let $\mathfrak M(V)$ denote the expected value of the random variable $N(f_1,\ldots,f_n)$.
Next, for the calculation of $\mathfrak M(V)$, we will need the notion of a Banach set on $X$,
as well as the notion of a volume of the Banach set.
\begin{definition}\label{dfBanach}
Let $T^*X$ be a cotangent bundle of $X$,
and ${\mathcal E}(x)$ be a convex centrally symmetric compact set in the cotangent space $T^*_xX$ of $X$ at the point $x$.
The collection ${\mathcal E} = \{{\mathcal E}(x)\subset T^*_xX \, \vert \, x \in X \}$
is called a \emph{Banach set} in $X$.
\end{definition}
\begin{definition}\label{dfBanachVol}
Consider the domain
$\bigcup _{x\in X}{\mathcal E }(x) \subset T^*(X).$
It's volume relative
to the standard symplectic structure in $T^*X$
is called the \emph{volume of Banach set} and is denoted by $\vol(\mathcal E)$.
\end{definition}
\noindent
For $x\in X$ let's define the linear functional $\theta(x)$ on $V$ as $\theta(x)(f)=f(x)$.
Next we assume that $\forall x\in X,\:\exists f\in V\colon\, f(x)\ne0$.
That's why the set $\theta(X)$ in $V^*$ does not contain $0$.
\begin{definition}\label{dfBanachOf}
Let's define the mapping $\Theta\colon X\to V^*$ as $\Theta(x)=\theta(x)/\sqrt{\langle\theta(x),\theta(x)\rangle_*}$,
where $\langle *,*\rangle_*$ is the scalar product in the space $V^*$ associated with the scalar product $\langle *,*\rangle$ in $V$.
Let $d\Theta_x\colon T_xX\to V^*$ be a differential of $\Theta$ at the point $x$.
Denote by $d^*\Theta_x\colon V\to T^*_xX$ the adjoint linear operator, and define the Banach set $\mathcal E_V$
by
$\mathcal E_V(x)=d^*\Theta_x(B)$,
where $B$ is the unit ball in $V$ centered at the origin.
The compact set $\mathcal E_V(x)$ is an ellipsoid.
We note that in a more general context discussed in [Ka1], arbitrary Banach sets on $X$ can arise.
\end{definition}
\begin{theorem}\label{thmAK1}
$\mathfrak M(V)=n!/(2\pi)^n\:\vol(\mathcal E_V)$
\end{theorem}
\begin{example}
Let $X$ be the unit circle $S^1$,
$V_m$ the space of trigonometric polynomials
$f(\theta)=\sum_{k\leq m} a_k\cos(k \theta)+b_k\sin(k \theta)$
 of degree $m$.
 Then (see \cite{ADG})
$$\mathfrak M(V_m)=\sqrt{\frac{m(m+1)}{3}}$$
For trigonometric polynomials in many variables see \cite{K22}.
\end{example}
Now let's state a similar theorem in the case where we consider $n$ spaces $V_i$ and equations $f_1=\ldots=f_n=0$, where $f_i\in V_i$.
For this, we will need the concept of the mixed volume of Banach sets.
Using Minkowski sum and homotheties, we
can form linear combinations of convex sets with non-negative coefficients.
The linear combination of Banach sets is defined by
$$(\sum _i \lambda _i {\mathcal E}_i)(x) = \sum _i \lambda _i {\mathcal E}_i(x).$$
For
$n$ Banach sets ${\mathcal E}_1, \ldots, {\mathcal E}_n$ the volume of $\lambda _1 {\mathcal E}_1+\ldots+\lambda _n{\mathcal E}_n$ is a homogeneous polynomial of degree $n$ in $\lambda _1, \ldots, \lambda _n$.
Its coefficient at $\lambda _1\cdot\ldots\cdot \lambda _n$ divided by $n!$ is called the mixed volume of Banach sets ${\mathcal E}_1, \ldots, {\mathcal E}_n$.
The mixed volume of Banach sets ${\mathcal E}_1, \ldots, {\mathcal E}_n$
is denoted by $\vol(\mathcal E_1,\ldots,\mathcal E_n)$.
\begin{theorem}\label{thmAK2}
Let $\mathfrak M(V_1,\ldots,V_n)$ denote the expectation of the random variable $N(f_1,\ldots,f_n)$.
Then it holds that
$$\mathfrak M(V_1,\ldots,V_n)=\frac{n!}{(2\pi)^n}\vol(\mathcal E_{V_1},\ldots,\mathcal E_{V_n})$$
\end{theorem}
\section*{The ring of Banach sets}
Next we need a concept of the ring of Banach sets.
It arises as an analogue of the well-known concept of a ring of convex bodies,
first defined in \cite{MCM}.
There are several different versions of this concept.
Here we construct an analogue of the definition from \cite{EKK}.
We call the formal difference $\mathcal E - \mathcal B$ of Banach sets the \emph{virtual Banach set}.
Virtual Banach sets form a vector space,
where multiplication by negative numbers is defined by $(-1)\cdot(\mathcal E - \mathcal B)=\mathcal B -\mathcal E$.

The following notations are used below

\noindent
$\bullet$ $S=\bigoplus_{0\leq i}S_i$ -- the graded symmetric algebra of the space of virtual Banach sets on the manifold $X$

\noindent
$\bullet$ $I$ -- the linear functional on the space $S$
defined by\ 1) $I_{S_k}=0$ for $k\ne n$, and\ 2) $I(\mathcal E_1\cdot\ldots\cdot\mathcal E_n)=\vol(\mathcal E_1,\ldots,\mathcal E_n)$

\noindent
$\bullet$ $L(x,y)=I(x\cdot y)$ --
the symmetric bilinear form on the vector space $S$

\noindent
$\bullet$ $J$ -- the kernel of the form $L$.
\begin{lemma}
$J$ is a homogeneous ideal of the graded ring $S$.
\end{lemma}
We will call ring $\mathfrak S=S/J$ the \emph{ring of virtual Banach sets}.
\begin{corollary}
The following statements hold:

\noindent
 {\rm(i)} $\mathfrak S_0=\R$

\noindent
 {\rm(ii)} $\dim \mathfrak S_n=1$

\noindent
 {\rm(iii)} The graded ring $\mathfrak S$ is generated by elements of degree $1$

\noindent
 {\rm(iv)} The mappings $\mathfrak S_p\times \mathfrak S_{n-p}\to\R,$ defined as $(\eta,\xi)\mapsto L(\eta,\xi)$, are non-degenerate pairings.
\end{corollary}

Next for $n$ virtual Banach sets $\mathcal B_1,\ldots,\mathcal B_n$,
we use the notation
$$
\vol(\mathcal B_1\cdot\ldots\cdot\mathcal B_n)
=I(\mathcal B_1\cdot\ldots\cdot\mathcal B_n)=\vol(\mathcal B_1,\ldots,\mathcal B_n)
$$
\section*{Zeros of random sections}
Transitioning to zeros of random sections of vector bundles, without formulating precise theorems,
we will first briefly describe the situation in the case
when considering zeros of sections of an $n$-dimensional vector bundle $\mathcal F$ on $X$.
Just as in the case of functions we consider
a finite dimensional space $V$ of smooth sections of $\mathcal F$.
Here we denote by $\mathfrak M(V;U)$ the expectation of the number of zeros of random section $s\in V$ contained in the open set $U\subset X$.
By ${\rm Res}_U\mathcal B$ we denote the constraint of $\mathcal B\in \mathfrak S$
on the subvariety $U$.
\begin{theorem}\label{thmSections}
There exists the unique element $\mathcal B\in \mathfrak S_n$,
such that for any $U\subset X$
$$\mathfrak M(V;U)=\frac{n!}{(2\pi)^n}\vol({\rm Res}_U\mathcal B)$$
\end{theorem}
Further results can be approximately described as follows.
We associate to an element $\mathfrak{s}$ of degree $k$ of the ring of Banach sets $\mathfrak{S}$ a certain $k$-density $d_k(\mathfrak s)$ on $X$
and interpret the ring $\mathfrak{S}$ as a ring of these densities.
Such densities serve as analogues of Chern forms, representing Chern classes of complex vector bundles,
and inherit some properties of Chern forms.

In conclusion, let us define the density $d_k(\mathfrak s)$.
An alternative construction of multiplication in the density ring is given in \cite{AK};
see also \cite{AKC}.
\begin{definition}\label{d_m} Let $H$ be the subspace of $T_xX$ generated by tangent vectors $\xi _1, \ldots, \xi_k$,
$H^\bot \subset T_x^*X$  the orthogonal complement to $H$, and $\pi _H: T_x^*X \to T_x^*X/H^\bot $
the projection map. The volume form on $T_x^*X/H^\bot $
is defined by $\omega(x)=\xi_1\wedge \ldots \wedge \xi_m$.
Let $\mathcal B_1,\ldots,\mathcal B_k$ be Banach sets on $X$.
Then
$d_k(\mathcal B_1\cdot\ldots\cdot\mathcal B_k)(\xi _1, \ldots, \xi _m)$ is the mixed
$k$-dimensional volume of convex $k$-dimensional sets $\pi_H\mathcal B_1(x), \ldots, \pi_H\mathcal B_k(x)$ in the sense of the volume form $\omega(x)$.
\end{definition}
\begin{lemma}
For different elements $\mathfrak s$ and $\mathfrak t$ of degree $k$ in the ring $\mathfrak S$,
the densities $d_k(\mathfrak s)$ and $d_k(\mathfrak t)$ are different.
\end{lemma}
Thus, the ring of virtual Banach sets $\mathfrak S$ can be considered as a certain ring of densities on the manifold $X$.

The following statement is an analogue of the BKK formula for Banach sets and for densities $d_i$.
\begin{theorem}\label{thm7, in1}
For any Banach sets $\mathcal B_1,\ldots,\mathcal B_k$  the equality
$$
d_1(\mathcal B_1)\cdot \ldots \cdot d_1(\mathcal B_k) =k!\,d_k(\mathcal B_1\cdot\ldots\cdot \mathcal B_k)
$$
holds.
\end{theorem}

\end{document}